\documentclass[a4paper]{article}
\usepackage{amsmath,amssymb,amsfonts}

\newtheorem{theorem}{Theorem}[section]
\newtheorem{definition}{Definition}[section]

\title{Infinite and natural numbers}

\author{Jailton C. Ferreira}

\date{ }
\begin{document}
\maketitle
\pagenumbering{arabic}
\begin{abstract}
The infinite numbers of the set $M$ of finite and infinite natural
numbers are defined starting from the sequence $0\Phi$, where $0$
is the first natural number, $\Phi$ a succession of symbols S and
$x$S is the successor of the natural number $x$. The concept of
limit of the natural number $n$, when $n$ tends to infinite, is
examined. Definitions and theorems about operations with elements
of $M$, equivalence and equality of natural numbers, distance
between elements of $M$ and the order of the elements are
presented.
\end{abstract}

\section{Introduction}

\hspace{22pt} Let us examine the natural numbers starting from the
$apeiron$ concept. In classical Greece $peiron/peras$ is what
delimits a thing from another (a finite boundary, a contour).
$Apeiron$ is what is limitless or endless. In Homer $apeiron$ is
endless and without internal distinctions and for Anaximander
$apeiron$ is the principle from which the $kosmoi$ are derived. If
we consider the conception, appeared in the VI B.C. century among
the pythagoreans, of the number as the essence of all the things,
we have that $apeiron$ is what does not have number in the
pythagoric sense, what does not have measure nor limit, the
immensurable, the formless, what does not have essence.

\hspace{22pt} Since $apeiron$ is infinite and formless, the
ancient Greeks considered $infinite$ as being imperfect and
incomplete. This meaning coexisted and prevailed on the idea of
infinite as perfection and completeness beyond any limit.

\hspace{22pt} The crucial point in Aristotle's conception of
infinite is that there is no transition from $kosmos$ to
$apeiron$, that infinite is not the $apeiron$ of his predecessors,
that it is not something formless and so involving that nothing
exists beyond $apeiron$ that is not itself. For Aristotle infinite
means that whatever is beyond does not distinguish from whatever
is on this side in terms of order. Aristotle allows infinite to
exist just in a restricted sense \cite{Aristotle}:

\begin{quote}
For generally the infinite has this mode of existence: one thing
is always being taken after another, and each thing that is always
finite, but always different.
\end{quote}

\hspace{22pt} The succession of finite natural numbers is endless,
that is, infinite. Considering that there is for all set the
cardinal number of the set, there is no way the set of the natural
numbers be a set of finite numbers and the number of elements of
the set be a natural number. Which the cardinality of the set of
digits that the natural numbers possess? In \cite{Ferreira1} we
showed that there are natural numbers with an endless number of
digits; we considered that the limit of a sequence $n$ of zeros,
when $n$ tends to infinite, is an infinite sequence of zeros.
Representing by $2^n$ the 1 preceded of $n$ zeros, in the binary
system, we considered

\begin{equation}\label{um-1}
\lim_{n \rightarrow \infty}2^{n}>j \qquad \textrm{for all}
\hspace{11pt} j\in N
\end{equation}

where $N$ is the set of the finite natural numbers.

\hspace{22pt} In the section \textbf{2} we defined the infinite
numbers of the set $M$ of finite and infinite natural numbers
starting from the sequence $0\Phi$, where 0 is the first natural
number, $\Phi$ a succession of symbols S and $x$S is the successor
of the natural number $x$. In the section \textbf{3} the limit
concept is extended and \eqref{um-1} explained.  Theorems about
operations with elements of $M$ are presented in the section
\textbf{4}, this section is contained in ~\cite{Ferreira2}. The
sections \textbf{5} to \textbf{7} treat of equivalence and
equality of numbers, distance between elements of $M$ and the
order of the elements. The section \textbf{8} contains a comment.

\hspace{22pt} As used in this article, $apeiron$ is neither the
chaos nor the indefinable, but it is what is indefinite and is, at
least partially, definable. $Apeiron$ is also associated with
ignorance or uncertainty about itself.

\section{The set M of natural numbers}

\hspace{22pt} Be the axioms
\begin{description}
\item[A1.] 0 is a natural number.
\item[A2.] For every natural number $x$ another number natural $x$S,
denominated successor of $x$, exists.
\end{description}

let us define the sequence, denoted by $y$, constituted of 0
followed by S's satisfying the following condition
\begin{gather}\label{dois-1}
\begin{tabular}{c c c c c c c}
$y$ &   & 0 & S & S & S & \ldots \\
  &   &   & $\updownarrow$ & $\updownarrow$ & $\updownarrow$ &  \\
  &   &   & 1 & 2 & 3 & \ldots \\
\end{tabular}
\end{gather}

where the numbers put into one-to-one correspondence with symbols
of $y$ belong to the set of the finite natural numbers, denoted by
$N$.

\hspace{22pt} Be the sequence $0\Delta$, such that $\Delta$ is a
sequence of S's and the set formed by S's of $\Delta$ is not
equivalent to any proper subset of itself. Let us notice that the
sequence $\Delta$ is indefinite, but $0\Delta$ is not formless:
there is in $0\Delta$ a certain order shared with finite natural
numbers. The sequence 0SSSS is a $0\Delta$ instance, it contains a
defined sequence, more ordinate, and it represents a natural
number.

\hspace{22pt} Let us consider the sequence such $0\Phi$ that
$\Phi$ is an endless sequence of S's. The sequence $0\Phi$ is more
indefinite than $y$, because in $0\Phi$ it is not made the
one-to-one correspondence between all S's of $\Phi$ and the
elements of $N$. In the same way that numbers belonging to $N$ are
instances of $0\Delta$, there are natural numbers no belonging to
$N$ that are $0\Phi$ instances.

\hspace{22pt} Be the sequence $0\Phi$ written as
\begin{gather}\label{dois-2}
\text{0$\mathfrak{SSS}$\ldots$\mathfrak{SSSS}$SSSSSSSS\ldots}
\end{gather}

where the symbols S indicated in the form $\mathfrak{S}$ are into
one-to-one correspondence with the elements of $N$ in the
following way
\begin{gather}\label{dois-3}
\begin{tabular}{c c c c c c c c c}
$\mathfrak{S}$ & $\mathfrak{S}$ & $\mathfrak{S}$ & $\mathfrak{S}$ & \ldots & $\mathfrak{S}$ & $\mathfrak{S}$ & $\mathfrak{S}$ & $\mathfrak{S}$ \\
$\updownarrow$ & $\updownarrow$ & $\updownarrow$ & $\updownarrow$ &   & $\updownarrow$ & $\updownarrow$ & $\updownarrow$ & $\updownarrow$\\
0 & 2 & 4 & 6 & \ldots & 7 & 5 & 3 & 1 \\
\end{tabular}
\end{gather}

Let us notice that the sequence
0$\mathfrak{SSS}$\ldots$\mathfrak{SSSS}$ is an instance of
0$\Phi$.

\begin{definition}
The sequence $\mathrm{0}\mathfrak{SSS}$\ldots$\mathfrak{SSSS}$ of
(\ref{dois-2}) is an infinite natural number.
\end{definition}

\hspace{22pt} Let us denote by $M$ the set of all the finite and
infinite natural numbers.

\begin{definition}

A1. $0$ is a natural number.

A2. For every natural number $x$ another number natural
$x\mathrm{S}$, denominated successor of $x$, exists.

A3. $0 \ne x\mathrm{S}$ for every natural number $x$.

A4. If $x\mathrm{S}$ = $y\mathrm{S}$, then $x$ = $y$.

A5. If a set of natural numbers, denoted by $M$, contains 0 and
also the successor of any natural number in $M$, then all the
natural numbers are in $M$.
\end{definition}

The Definition 2.2 generalizes the Peano's axioms substituting $N$
by $M$.

\begin{definition}
Be $K$ a number, greater than zero, belonging to the set of
natural numbers $M$
\begin{gather}\label{dois-4}
\text{$M=\{0,1,2,3,4, \ldots,
K_{-3},K_{-2},K_{-1},K,K_{1},K_{2},K_{3}, \dots\}$ }
\end{gather}
where $K_n$ denotes the $n$-th successor of $K$ and $K_{-n}$
denotes the $n$-th predecessor of $K$. If $K$ is such that the set
\begin{gather}\label{dois-5}
\text{$\{K,K_{-1},K_{-2},K_{-3}, \ldots, 3,2,1,0\}$ }
\end{gather}

is equivalent to $N$, that is,
\begin{gather}\label{dois-6}
\begin{tabular}{c c c c c c c c c}
$\{$ & $K,$ & $0,$ & $K_{-1},$ & $1,$ & $K_{-2},$ & $2,$ & $K_{-3},$ & $\ldots\}$\\
  & $\updownarrow$ & $\updownarrow$ & $\updownarrow$ & $\updownarrow$ & $\updownarrow$ & $\updownarrow$ & $\updownarrow$ &  \\
$\{$ & $0,$ & $1,$ & $2,$ & $3,$ & $4,$ & $5,$ & $6,$ & $\ldots\}$\\
\end{tabular}
\end{gather}

then $K$ is an infinite number.
\end{definition}

$K$ and the sequence 0$\mathfrak{SSS}\ldots\mathfrak{SSSS}$ of
($\ref{dois-2}$) are two general forms of representing an infinite
natural number.

\hspace{22pt} Let us consider the following instance of $0\Phi$:
\begin{gather}\label{dois-7}
\text{0$\mathfrak{S}\mathrm{SS}\mathfrak{SSSS}\mathrm{SSSSSSSS}\mathfrak{SSSSSSSSSSSSSSSS}\ldots$}
\end{gather}

This instance can be reordained in the form of the column on the
right of the list
\begin{gather}\label{dois-8}
\begin{tabular}{r l}
0 & $0\mathfrak{S}$\\
1 & $\hspace{5pt}\mathrm{SS}$\\
2 & $\hspace{5pt}\mathfrak{SSSS}$\\
3 & $\hspace{5pt}\mathrm{SSSSSSSS}$\\
4 & $\hspace{5pt}\mathfrak{SSSSSSSSSSSSSSSS}$\\
$\ldots$ & $\hspace{5pt}\ldots$\\
\end{tabular}
\end{gather}

In (\ref{dois-8}) the elements of the column on the left are the
finite natural numbers and the $i$+1-th sequence of S's has the
double of symbols that the $i$-th sequence of S's. The right
column of (\ref{dois-8}) can be represented by the sequence of
symbols
\begin{gather}\label{dois-9}
\text{\ldots11111111111}
\end{gather}

where the first 1 from the right to the left corresponds to S of
the row 1 of (\ref{dois-8}), the second 1 from the right to the
left corresponds to S's of the row 2 of (\ref{dois-8}), the $n$-th
1 from the right to the left corresponds to S's of the row $n$ of
(\ref{dois-8}), where $n \in N$, and so forth.

\hspace{22pt}The form (\ref{dois-9}) is an instance of
\begin{gather}\label{dois-10}
\text{$\ldots f_5 f_4 f_2 f_2 f_1 $}
\end{gather}

where $f_i \in \{0, 1\}$ and $i \in N$.

To a representation (\ref{dois-10}) with $f_i = 0$ corresponds an
instance of 0$\Phi$, such that the list in the form (\ref{dois-9})
possesses the $i$-th row empty of symbols S. The instances of
(\ref{dois-10}) in which for any $f_i$ exists some $f_j$ such that
$j > i$ and $f_j = 1$ are instances of $0\Phi$. An instance of
(\ref{dois-10}) in the which exists some $i \in N$ such that $f_j
= 0$ for $j > i$, is a simplified form of representing
\begin{gather}\label{dois-11}
\text{$f_i \times 2^{i-1} + \ldots + f_4 \times 2^3 +  f_3 \times
2^2 +f_2 \times 2^1 +f_1 \times 2^0 $}
\end{gather}

Let us also notice that the form (\ref{dois-10}) does not
represent all instances of $0\Phi$, as shows the following
theorem:
\begin{theorem}
 The successor of the natural number represented in the form
(\ref{dois-10}) such that $f_i = 1$ for every $i \in N$, is not
representable in the form (\ref{dois-10}).
\end{theorem}

\textit{Proof:}

\hspace{22pt} Be the natural number represented in the form
(\ref{dois-10}) such that $f_i=1$ for every $i\in N$, we
denominate this number $w$. Since every $f_i$ is equal to 1, it
does not exist $f_i$ to be changed from 0 to 1 in the
representation of $w$S in the form (\ref{dois-10}). Therefore, for
$w$S no $f_i$ is equal to 0. If $w$S is such that $f_i=1$ for
every $i\in N$, then $w$S$=w$. Since a natural number cannot be
the immediate successor of itself, we concluded that $w$S is not
representable in the form
(\ref{dois-10}).\\
The successors of $w$ are indicated by
\begin{gather}\label{dois-12}
\text{$w_1, w_2 , w_3 , \ldots$}
\end{gather}

therefore
\begin{gather}\label{dois-13}
\text{$M=\{0,1,2,\ldots , w,w_1, w_2 , w_3 , \ldots \}$}
\end{gather}

\section{Concepts of limit and M}

\hspace{22pt}Besides the limit definition by Heine and Weierstrass

\begin{quote}
If, given any $\epsilon$, there is a $\delta_0$ such that for
$0<\delta<\delta_0$ the difference $f(x_0 \pm \delta) - L$ is
smaller in absolute value than $\epsilon$, then $L$ is the limit
of $f(x)$ for $x = x_0$.
\end{quote}

the word $limit$ usually indicates (i) the existence of
boundaries, (ii) exclusion which either are not passed over or
cannot or may not be (iii) part or extreme point.

\hspace{22pt} Let us adopt $limit$ firstly indicating point not to
be passed. Be it denoted by
\begin{gather}\label{tres-1}
\text{$1 \ldots _{n} 11$}
\end{gather}

the sequence containing $n$ symbols 1 with $i \in N$ and by $L$
the set constituted by the elements of $M$ represented in the form
(\ref{dois-10}). What is the limit of (\ref{tres-1}) when the
number $n$ of symbols 1, being $n$ belonging to $L$, does tend to
infinite? This limit is $ \ldots f_5 f_4 f_3 f_2 f_1$ with $f_i =
1$ for every $i \in N$. Therefore
\begin{equation}\label{tres-2}
\lim_{(n \in L, n \rightarrow \infty )}1 \ldots _{n} 11 = w
\end{equation}

The leminiscate $\infty$ just means $limitless$. Let us also
notice that the limit does not exist
\begin{equation}\label{tres-3}
\lim_{(n \in N, n \rightarrow \infty )}n
\end{equation}
and that

\begin{equation}\label{tres-4}
\lim_{(n \in L, n \rightarrow \infty )}2^{n} = \lim_{(n \in L, n
\rightarrow \infty )}[(1 \ldots _{n} 11) + 1] =[\lim_{(n \in L, n
\rightarrow \infty )}(1 \ldots _{n} 11)] + 1 = w_1
\end{equation}

\hspace{22pt} Let us now examine $limit$ indicating extreme part.
What is the limit of $n$ belonging to $M$ when $n$ tends to
infinite? Something that is beyond $M$ is not a natural number, it
is what there is in $0\Phi$ that cannot be instantiated as natural
number. In certain a way
\begin{equation}\label{tres-5}
\lim_{xtr(n \in M, n \rightarrow \infty )}n
\end{equation}

where $xtr$ denotes $limit$ as extreme part, is an instance of the
concept $apeiron$. Let us notice that

\begin{equation}\label{tres-6}
\lim_{xtr(n \in N, n \rightarrow \infty )}n = K
\end{equation}

\section{Operations with $K$}

\begin{definition}
Be $0X$ such that $X$ is a sequence of S's. The natural number $n$
is the cardinal number of a set $A$ when there is one-to-one
correspondence between the elements of $A$ and the symbols S of
the sequence $0X$ that represents the number $n$.
\end{definition}

The cardinal number of $A$ is denoted by $\lvert A \lvert$. In the
case of the set of natural numbers $N$, we have

\begin{gather}\label{quatro-1}
\begin{tabular}{l c r c c c c l}
$N$ & & $\{ 0,$ & $1,$ & $2,$ & $3,$ & $4,$ & $\ldots \}$ \\
  &   & $\updownarrow$ & $\updownarrow$ & $\updownarrow$ & $\updownarrow$ & $\updownarrow$ &  \\
$\lvert N \lvert$ & 0 & S & S & S & S & S & $\ldots$ \\
\end{tabular}
\end{gather}

Let us notice that
\begin{gather}\label{quatro-2}
\text{$K=\lvert N \lvert$}
\end{gather}

\begin{definition}
The every pair of natural numbers $x,y$, we can attribute a
natural number, denominated $x+y$, such that
\begin{enumerate}
\item $x$ is the cardinal number of a set $A$
\item $y$ is the cardinal number of a set $B$
\item $A \cap B = \varnothing$
\item $x+y$ is the cardinal number of the set $A \cup
B$
\end{enumerate}
$x+y$ is denominated the sum of $x$ and $y$, or the number
obtained by the addition of $y$ to $x$.
\end{definition}

\begin{theorem}
The sum $K+n$  is equal to $K$.
\end{theorem}

\textit{Proof:}\\
\hspace{22pt}Let the disjoint sets
\begin{gather}\label{quatro-3}
\text{$A= \{ a_1, a_2, a_3, \ldots , a_n\}$}
\end{gather}
where $n \in N$, and
\begin{gather}\label{quatro-4}
\text{$B= \{ K, K_{-1}, K_{-2}, K_{-3}, \ldots , 3,2,1,0\}$}
\end{gather}

From (\ref{quatro-2}), (\ref{quatro-3}) and (\ref{quatro-4}) we
have
\begin{gather}\label{quatro-5}
\text{$\lvert A \lvert = n$}
\end{gather}
and
\begin{gather}\label{quatro-6}
\text{$\lvert B \lvert = \lvert N \lvert = K$}
\end{gather}

We can establish the following one-to-one correspondence of $A
\cup B$ with $N$
\begin{gather}\label{quatro-7}
\begin{tabular}{c c c c c c c c c c}
$\{ a_1 ,$ & $a_2 ,$& $a_3 ,$ & $\ldots ,$ & $a_n ,$ & $K,$ & $0,$ & $K_{-1},$ & $1,$ & $\ldots \}$ \\
$\updownarrow$ & $\updownarrow$ & $\updownarrow$ &   & $\updownarrow$ & $\updownarrow$ & $\updownarrow$ & $\updownarrow$ &  $\updownarrow$ & \\
$\{ 0,$ & $1,$ & $2,$ & $\ldots ,$ & $n+1$ & $n+2$ & $n+3$ & $n+4$ & $n+5$ & $\ldots \}$ \\
\end{tabular}
\end{gather}

therefore
\begin{gather}\label{quatro-8}
\text{$\lvert A \cup B\lvert = \lvert N \lvert$}
\end{gather}

Since $A$ and $B$ are disjoint sets, we can affirm
\begin{gather}\label{quatro-9}
\text{$\lvert A \lvert + \lvert B \lvert = \lvert A \cup B\lvert$}
\end{gather}

Substituting (\ref{quatro-5}), (\ref{quatro-6}) and
(\ref{quatro-8}) in (\ref{quatro-9}), we conclude
\begin{gather}\label{quatro-10}
\text{$K+n=K$}
\end{gather}

\begin{theorem}
The difference $K-n$ is equal to $K$.
\end{theorem}

\textit{Proof:}

\hspace{22pt} Let
\begin{gather}\label{quatro-11}
\text{$A= \{ K, K_{-1}, K_{-2}, K_{-3}, \ldots , 3,2,1,0\}$}
\end{gather}
and
\begin{gather}\label{quatro-12}
\text{$B= \{ K, K_{-1}, K_{-2}, K_{-3}, \ldots , K_{-n} \}$}
\end{gather}

where $n \in N$. Let us consider
\begin{gather}\label{quatro-13}
\text{$\lvert A \lvert = K$}
\end{gather}
\begin{gather}\label{quatro-14}
\text{$\lvert B \lvert = n$}
\end{gather}
and
\begin{gather}\label{quatro-15}
\text{$\lvert A-B \lvert = \lvert A \lvert - \lvert B \lvert$}
\end{gather}

We can establish the following one-to-one correspondence of $A-B$
with $N$:
\begin{gather}\label{quatro-16}
\begin{tabular}{c c c c c c c c c}
$\{ K_{-(n-1)} ,$ & $0 ,$& $K_{-(n-2)} ,$ & $1 ,$ & $K_{-(n-3)} ,$ & $2,$ & $K_{-(n-4)} ,$ & $3,$ & $\ldots \}$ \\
$\updownarrow$ & $\updownarrow$ & $\updownarrow$ & $\updownarrow$ & $\updownarrow$ & $\updownarrow$ & $\updownarrow$ & $\updownarrow$ &  \\
$\{ 0,$ & $1,$ & $2,$ & $3,$ & $4,$ & $5,$ & $6,$ & $7,$ & $\ldots \}$ \\
\end{tabular}
\end{gather}

that is
\begin{gather}\label{quatro-17}
\text{$\lvert A-B \lvert = \lvert N \lvert = K$}
\end{gather}

Substituting (\ref{quatro-13}), (\ref{quatro-14}) and
(\ref{quatro-17}) in (\ref{quatro-15}), we conclude
\begin{gather}\label{quatro-18}
\text{$K-n=K$}
\end{gather}

\begin{theorem}
The sum $K+K$ is equal to $K$.
\end{theorem}

\textit{Proof:}

\hspace{22pt} Let the disjoint sets
\begin{gather}\label{quatro-19}
\text{$A= \{ K, K_{-1}, K_{-2}, K_{-3}, \ldots , 3,2,1,0\}$}
\end{gather}
and
\begin{gather}\label{quatro-20}
\text{$B= \{ b_1, b_2, b_3, \ldots , b_K\}$}
\end{gather}

Let us note
\begin{gather}\label{quatro-21}
\text{$\lvert A \lvert = K$}
\end{gather}
and
\begin{gather}\label{quatro-22}
\text{$\lvert B \lvert = K$}
\end{gather}

We can establish the following one-to-one correspondence of $A
\cup B$ with $N$:
\begin{gather}\label{quatro-23}
\begin{tabular}{c c c c c c c c c}
$\{ K,$ & $0 ,$& $b_K ,$ & $b_1 ,$ & $K_{-1} ,$ & $1,$ & $b_{K_{-1}} ,$ & $b_2,$ & $\ldots \}$ \\
$\updownarrow$ & $\updownarrow$ & $\updownarrow$ & $\updownarrow$ & $\updownarrow$ & $\updownarrow$ & $\updownarrow$ & $\updownarrow$ &  \\
$\{ 0,$ & $1,$ & $2,$ & $3,$ & $4,$ & $5,$ & $6,$ & $7,$ & $\ldots \}$ \\
\end{tabular}
\end{gather}
therefore
\begin{gather}\label{quatro-24}
\text{$\lvert A \cup B \lvert = \lvert N \lvert$}
\end{gather}

Since $A$ and $B$ are disjoint sets, we have
\begin{gather}\label{quatro-25}
\text{$\lvert A \lvert + \lvert B \lvert = \lvert A \cup B
\lvert$}
\end{gather}

Substituting (\ref{quatro-21}), (\ref{quatro-22}),
(\ref{quatro-24}) and (\ref{quatro-2}) in (\ref{quatro-25}), we
conclude
\begin{gather}\label{quatro-26}
\text{$K+K=K$}
\end{gather}

\begin{definition}
What $\kappa$ represents can be either number $K$ or a finite
natural number.
\end{definition}

\begin{theorem}
The difference $K-K$ is equal to $\kappa$.
\end{theorem}

\textit{Proof:}

\hspace{22pt} Let the sets
\begin{gather}\label{quatro-27}
\text{$A= \{ a_0, a_1, a_2, a_3, \ldots , a_K\}$}
\end{gather}
\begin{gather}\label{quatro-28}
\text{$B= \{ b_0, b_1, b_2, b_3, \ldots , b_K\}$}
\end{gather}
where
\begin{gather}\label{quatro-29}
\text{$B \subset A$}
\end{gather}
and
\begin{gather}\label{quatro-30}
\text{$C=A-B$}
\end{gather}

If the elements of $B$ can be put in the sequence of elements of
$A$, so that no element of $A$ belonging to $B$ immediately
succeeds or immediately precedes another element of $A$ belonging
to $B$; then each set of elements of $A$ belonging to $C$ between
two elements of $A$ belonging to $B$ can be put into one-to-one
correspondence with a natural number and
\begin{gather}\label{quatro-31}
\text{$\lvert C \lvert = K$}
\end{gather}

For instance, for
\begin{gather}\label{quatro-32}
\text{$A=\{ c_0, c_5, b_0, c_1, b_1, c_2, c_9, b_2, c_{20}, \ldots
, c_4, b_3, \ldots \} $}
\end{gather}

where $c_i$ belongs to $C$, we have
\begin{gather}\label{quatro-33}
\begin{tabular}{c c c c c c}
$\{ $ & $c_0,c_5,$& $c_1,$ & $c_2,c_9,$ & $c_{20}, \ldots ,c_4,$ & $\ldots \}$ \\
  & $\updownarrow$ & $\updownarrow$ & $\updownarrow$ & $\updownarrow$ &  \\
$\{$ & $0,$ & $1,$ & $2,$ & $3,$ & $\ldots \}$ \\
\end{tabular}
\end{gather}

Considering
\begin{gather}\label{quatro-34}
\text{$\lvert A-B \lvert = \lvert A \lvert - \lvert B \lvert$}
\end{gather}
\begin{gather}\label{quatro-35}
\text{$\lvert A \lvert = K$}
\end{gather}
\begin{gather}\label{quatro-36}
\text{$\lvert B \lvert = K$}
\end{gather}
we conclude that, after we substitute (\ref{quatro-30}),
(\ref{quatro-31}), (\ref{quatro-35}) and (\ref{quatro-36}) in
(\ref{quatro-34}),
\begin{gather}\label{quatro-37}
\text{$K-K=K$}
\end{gather}

\hspace{22pt}If the elements of $A$ can be put in the form
\begin{gather}\label{quatro-38}
\text{$A= \{ c_0, \ldots , c_n, b_0, b_1, b_2, b_3, \ldots \}$}
\end{gather}
with
\begin{gather}\label{quatro-39}
\text{$\lvert C \lvert = n$ and $n \in N$}
\end{gather}

then substituting (\ref{quatro-30}), (\ref{quatro-35}),
(\ref{quatro-36}) and (\ref{quatro-39}) in (\ref{quatro-34}), we
have
\begin{gather}\label{quatro-40}
\text{$K-K= n$}
\end{gather}

From (\ref{quatro-37}) and (\ref{quatro-40}) we have $K-K$ is
either $K$ or a finite number, that is,
\begin{gather}\label{quatro-41}
\text{$K-K= \kappa$}
\end{gather}

\begin{theorem}
The product $K \times n$ is equal to $K$.
\end{theorem}
\textit{Proof:}

\hspace{22pt} Let ${(K + )}_n$ be the sum of the $n$ portions
equal to $K$ and $n \in N$

(i) For $n = 1$, we have
\begin{gather}\label{quatro-42}
\text{$K \times 1 = {(K + )}_1 = K$}
\end{gather}

(ii) Let us assume
\begin{gather}\label{quatro-43}
\text{$K \times n = {(K + )}_n = K$}
\end{gather}

(iii)
\begin{gather}\label{quatro-44}
\text{$K \times (n+1) = {(K + )}_{n+1}$}
\end{gather}
\begin{gather}\label{quatro-45}
\text{$K \times (n+1) = {(K + )}_n + K$}
\end{gather}
\begin{gather}\label{quatro-46}
\text{$K \times (n+1) = K + K$}
\end{gather}
applying the Theorem 4-3 to (\ref{quatro-46}), we have
\begin{gather}\label{quatro-47}
\text{$K \times (n+1) = K$}
\end{gather}

(iv) If $K \times n=K$ holds for $n$, then $K×(n+1)=K$ holds

(v) Therefore, from steps (i) and (iv), we have $K \times n=K$ for
every $n$.

\begin{theorem}
The quotient $K \div n$ is equal to $K$.
\end{theorem}
\textit{Proof:}

\hspace{22pt} Let the sets
\begin{gather}\label{quatro-48}
\begin{tabular}{c c c c l}
$A=$ & $\{ (c_1, b_{1}),$ & $(c_2, b_{1}),$ & $(c_3, b_{1}),$ & $\ldots$ \\
  & $(c_1, b_{2}),$ & $(c_2, b_{2}),$ & $(c_3, b_{2}),$ & $\ldots$ \\
  & $\ldots$ & $\ldots$ & $\ldots$ & $\ldots$ \\
  & $(c_1, b_{n}),$ & $(c_2, b_{n}),$ & $(c_3, b_{n}),$ & $\ldots \}$ \\
\end{tabular}
\end{gather}

\begin{gather}\label{quatro-49}
\text{$B= \{ b_1, b_2, b_3, \ldots , b_n\}$}
\end{gather}
and
\begin{gather}\label{quatro-50}
\text{$C= \{ c_1, c_2, c_3, \ldots \}$}
\end{gather}
such that
\begin{gather}\label{quatro-51}
\text{$\lvert A \lvert = K$}
\end{gather}
and
\begin{gather}\label{quatro-52}
\text{$\lvert B \lvert = n > 0$ and $n \in N$}
\end{gather}

The quotient $K \div n$ is the number of elements per row of
(\ref{quatro-48}); considering that each row constitutes an
equivalent set to $A$, we have
\begin{gather}\label{quatro-53}
\text{$K \div n = K$}
\end{gather}

\begin{theorem}
The product $K \times K$ is equal to $K$.
\end{theorem}
\textit{Proof:}

\hspace{22pt} Let the sets
\begin{gather}\label{quatro-54}
\begin{tabular}{c c c c l}
$A=$ & $\{ (b_1, c_{1}),$ & $(b_1, c_{2}),$ & $(b_1, c_{3}),$ & $\ldots$ \\
  & $(b_2, c_{1}),$ & $(b_2, c_{2}),$ & $(b_2, c_{3}),$ & $\ldots$ \\
  & $(b_3, c_{1}),$ & $(b_3, c_{2}),$ & $(b_3, c_{3}),$ & $\ldots$ \\
  & $\ldots$ & $\ldots$ & $\ldots$ & $\ldots \}$ \\
\end{tabular}
\end{gather}
\begin{gather}\label{quatro-55}
\text{$B= \{ b_1, b_2, b_3, \ldots$}
\end{gather}
and
\begin{gather}\label{quatro-56}
\text{$C= \{ c_1, c_2, c_3, \ldots \}$}
\end{gather}
such that
\begin{gather}\label{quatro-57}
\text{$\lvert B \lvert = K$}
\end{gather}
and
\begin{gather}\label{quatro-58}
\text{$\lvert C \lvert = K$}
\end{gather}

Let us simplify the representation of the elements ($b_i,c_j$) to
($i, j$). Disposing the elements of $A$ as in the Table 1 and
applying the diagonalizing technique, we can establish the
one-to-one correspondence between the pairs ($i, j$) with $N$.
\begin{gather}\label{quatro-59}
\begin{tabular}{c c c c c c c c c c}
$A$ &  & $(1,1)$ & $(1,2)$ & $(2,1)$ & $(3,1)$ & $(2,2)$ & $(1,3)$ & $(1,4)$ & $\ldots$ \\
  &  & $\updownarrow$ & $\updownarrow$ & $\updownarrow$ & $\updownarrow$ & $\updownarrow$ & $\updownarrow$ & $\updownarrow$ &  \\
$N$ &  & 0 & 1 & 2 & 3 & 4 & 5 & 6 & $\ldots$ \\
\end{tabular}
\end{gather}

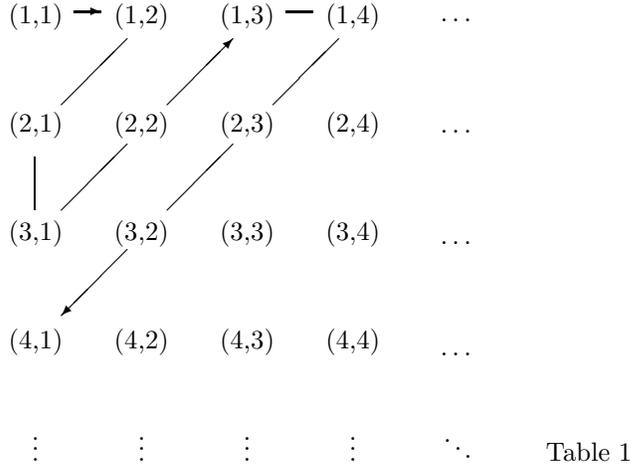
\begin{center}
\begin{picture}(200, 200)
\put(0,0){\shortstack{(1,1)\vspace{28pt} \\ (2,1)\vspace{28pt} \\ (3,1)\vspace{28pt} \\
(4,1)\vspace{22pt} \\ \vdots}}
\put(40,0){\shortstack{(1,2)\vspace{28pt} \\ (2,2)\vspace{28pt} \\ (3,2)\vspace{28pt} \\
(4,2)\vspace{22pt} \\ \vdots}}
\put(80,0){\shortstack{(1,3)\vspace{28pt} \\ (2,3)\vspace{28pt} \\ (3,3)\vspace{28pt} \\
(4,3)\vspace{22pt} \\ \vdots}}
\put(120,0){\shortstack{(1,4)\vspace{28pt} \\ (2,4)\vspace{28pt} \\ (3,4)\vspace{28pt} \\
(4,4)\vspace{22pt} \\ \vdots}}
\put(160,0){\shortstack{\ldots \vspace{38pt} \\ \ldots \vspace{38pt} \\ \ldots \vspace{38pt} \\
\ldots \vspace{22pt} \\ $\ddots$}}

\put (45,80) {\vector(-1,-1){25}} \put (60,95) {\line(1,1){25}}
\put (100,135) {\line(1,1){25}}

\put (20,95) {\line(1,1){25}} \put (60,135) {\vector(1,1){25}}

\put (20,135) {\line(1,1){25}}

\put (25,170) {\vector(1,0){10}}

\put (105,170) {\line(1,0){10}}

\put (10,95) {\line(0,1){20}}
\end{picture}
\nolinebreak \text{Table 1}
\end{center}

The number of rows of $A$ is $K$ and the number of columns is $K$;
therefore, the number of elements of $A$ is $K \times K$. Since
the number of elements of $A$ is  $\lvert N \lvert$, we conclude
that
\begin{gather}\label{quatro-60}
\text{$K \times K = K$}
\end{gather}

\begin{theorem}
The quotient $K \div K$ is equal to $\kappa$.
\end{theorem}
\textit{Proof:}

\hspace{22pt}   Let the sets
\begin{gather}\label{quatro-61}
\begin{tabular}{c c c c l}
$A=$ & $\{ (b_1, c_{1}),$ & $(b_1, c_{2}),$ & $(b_1, c_{3}),$ & $\ldots$ \\
  & $(b_2, c_{1}),$ & $(b_2, c_{2}),$ & $(b_2, c_{3}),$ & $\ldots$ \\
  & $(b_3, c_{1}),$ & $(b_3, c_{2}),$ & $(b_3, c_{3}),$ & $\ldots$ \\
  & $\ldots$ & $\ldots$ & $\ldots$ & $\ldots \}$ \\
\end{tabular}
\end{gather}
\begin{gather}\label{quatro-62}
\text{$B= \{ b_1, b_2, b_3, \ldots \}$}
\end{gather}
and
\begin{gather}\label{quatro-63}
\text{$C= \{ c_1, c_2, c_3, \ldots \}$}
\end{gather}
such that
\begin{gather}\label{quatro-64}
\text{$\lvert A \lvert = K$}
\end{gather}
\begin{gather}\label{quatro-65}
\text{$\lvert B \lvert = K$}
\end{gather}
and
\begin{gather}\label{quatro-66}
\text{$\lvert C \lvert = K$}
\end{gather}

The number of elements per row of (\ref{quatro-61}) is $K \div K$.
Considering that each row constitutes an set equivalent to $C$ and
(\ref{quatro-66}), we have
\begin{gather}\label{quatro-67}
\text{$K \div K = K$}
\end{gather}
If
\begin{gather}\label{quatro-68}
\text{$\lvert C \lvert = n$ and $n \in N$}
\end{gather}
we have
\begin{gather}\label{quatro-69}
\begin{tabular}{r c c c r}
$A=$ & $\{ (b_1, c_{1}),$ & $(b_1, c_{2}),$ & $\ldots ,$ & $(b_1, c_{n})$ \\
  & $(b_2, c_{1}),$ & $(b_2, c_{2}),$ & $\ldots ,$ & $(b_2, c_{n})$ \\
  & $(b_3, c_{1}),$ & $(b_3, c_{2}),$ & $\ldots ,$ & $(b_3, c_{n})$ \\
  & $\ldots$ & $\ldots$ & $\ldots$ & $\ldots \}$ \\
\end{tabular}
\end{gather}

The number of elements per row of (\ref{quatro-69}) is $K \div K$.
Considering that each row constitutes a set of cardinality $n$, we
have
\begin{gather}\label{quatro-70}
\text{$K \div K = n$}
\end{gather}

From (\ref{quatro-67}) and (\ref{quatro-70}) we have that $K \div
K$ can be either $K$ or a finite number, that is,
\begin{gather}\label{quatro-71}
\text{$K \div K = \kappa$}
\end{gather}

\begin{theorem}
$K + \kappa = K$
\end{theorem}
\textit{Proof:}

\hspace{22pt} If $\kappa$ is equal to $K$, then $K+K=K$ (Theorem
4.3); if $\kappa$ is equal to $n$ and $n \in N$, then $K+n=K$
(Theorem 4.1). Therefore, $K + \kappa = K$.

\begin{theorem}
$K - \kappa = \kappa$
\end{theorem}
\textit{Proof:}

\hspace{22pt} If $\kappa$ is equal to $K$, then $K-K=\kappa$
(Theorem 4.4); if $\kappa$ is equal to $n$ and $n \in N$, then
$K-n=K$ (Theorem 4.2). Therefore, $K - \kappa = \kappa$.

\begin{theorem}
$K \times \kappa = K$
\end{theorem}
\textit{Proof:}

\hspace{22pt} If $\kappa$ is equal to $K$, then $K \times K = K$
(Theorem 4.7); if $\kappa$ is equal to $n$ and $n \in N$, then $K
\times n = K$ (Theorem 4.5). Therefore, $K \times \kappa = K$.

\begin{theorem}
$K \div \kappa = \kappa$
\end{theorem}
\textit{Proof:}

\hspace{22pt} If $\kappa$ is equal to $K$, then $K \div K =
\kappa$ (Theorem 4.8); if $\kappa$ is equal to $n$ and $n \in N$,
then $K \div n = K$ (Theorem 4.6). Therefore, $K \div \kappa =
\kappa$

\begin{theorem}
Be $z \in (M-N)$ and identified by a structural-descriptive name
in the form (\ref{dois-10}). The sum $z + K$ is equal to $K$.
\end{theorem}
\textit{Proof:}

\hspace{22pt} $z$ represents a given number of $M$, while $K$
represents any infinite number of $M$. From the Theorem 4.3 we
have
\begin{gather}\label{quatro-72}
\text{$K + K = K$}
\end{gather}

instantiating the first portion of (\ref{quatro-72}) to $z$, we
concluded
\begin{gather}\label{quatro-73}
\text{$z + K = K$}
\end{gather}

\begin{theorem}
Be $x$ and $y$ belonging to $M-N$ and identified by
structural-descriptive names in the form (\ref{dois-10}). The sum
$x + y$ is equal to $K$.
\end{theorem}
\textit{Proof:}

\hspace{22pt} $x$ and $y$ represent two given numbers of $M$,
while $K$ represents any infinite number of $M$. Of the Theorem
4.3 we have
\begin{gather}\label{quatro-74}
\text{$K + K = K$}
\end{gather}

instantiating the first portion of (\ref{quatro-74}) for $x$ and
second parcel to $y$, we concluded
\begin{gather}\label{quatro-75}
\text{$x + y = K$}
\end{gather}

\section{Equivalence and equality of numbers}

\begin{definition}
Be (i) $n$ and $m$ belonging to $M$, (ii) $X_n$ and $X_m$
sequences of S's and (iii) the sequences 0$X_n$ and 0$X_m$
representations of $n$ and $m$ in the form $0X$, respectively. If
the set of symbols S of $X_n$ is a proper subset of the set of
symbols S of $X_m$, then $n <m$. If the set of symbols S of $X_m$
is a proper subset of the set of symbols S of $X_n$, then $n> m$.
If the set of symbols S of $X_n $ is equivalent to the set of
symbols S of $X_m$, then $n $ is equivalent the $m$.
\end{definition}

\begin{theorem}
If $m \in M$, $n \in M$ and $n$ is equivalent to $K$, then $n $ is
greater or equivalent to $m$.
\end{theorem}
\textit{Proof:}

\hspace{22pt} If $m \in N$, we can write
\begin{gather}\label{cinco-1}
\begin{tabular}{c c c c c c c c c c}
0 & $S$ & $S$ & $S$ & $\ldots$ &  $S$ & $\ldots$ & $S$ & $S$ & $S$ \\
  & $\updownarrow$ & $\updownarrow$ & $\updownarrow$ &   &  $\updownarrow$ &   & $\updownarrow$ & $\updownarrow$ & $\updownarrow$ \\
  & $\{ 0,$ & $1,$ & $2,$ & $\ldots ,$ &  $m,$ & $\ldots ,$ & $K_{-2}$ & $K_{-1}$ & $K \}$ \\
  & $\updownarrow$ & $\updownarrow$ & $\updownarrow$ &   &  $\updownarrow$ &   &   &   &   \\
  & $\{ 0,$ & $1,$ & $2,$ & $\ldots ,$ &  $m \}$ &  &  &  &  \\
  & $\updownarrow$ & $\updownarrow$ & $\updownarrow$ &   &  $\updownarrow$ &   &   &   &   \\
0 & $S$ & $S$ & $S$ & $\ldots$ &  $S$ &   &   &  &  \\
\end{tabular}
\end{gather}

The first and last rows of (\ref{cinco-1}) are $n$ and $m$
represented in the form $0X$. If $m \notin N$, then so much $m$ as
$n$ are represented in the form $0X$ by (\ref{cinco-2}).
\begin{gather}\label{cinco-2}
\begin{tabular}{c c c c c c c c}
0 & $S$ & $S$ & $S$ & $\ldots$ &  $S$ &  $S$ & $S$ \\
  & $\updownarrow$ & $\updownarrow$ & $\updownarrow$ &   & $\updownarrow$ & $\updownarrow$ & $\updownarrow$ \\
  & $\{ 0,$ & $1,$ & $2,$ & $\ldots ,$ &  $K_{-2}$ & $K_{-1}$ & $K \}$ \\
\end{tabular}
\end{gather}

From (\ref{cinco-1}), (\ref{cinco-2}) and Definition 5.1, we
concluded that $n$ is greater or equivalent to $m$.

\begin{definition}
Be $n$ and $m$ belonging to $N$. If $n$ is equivalent to $m$, then
$n$ is equal to $m$, that is, $n = m$.
\end{definition}

\begin{definition}
Be $n$ and $m$ belonging to $M$ and $n$ equivalent to $m$. If the
structural-descriptive names of $n$ and $m$ identify the same
number, then $n$ is equal to $m$, that is, $m = n$.
\end{definition}

\section{Distance between elements}

Let us denote by $P$ the set of the elements belonging to $M-N$
represented in the form (\ref{dois-10}).

\begin{definition}
Be $x$ and $y$ belonging to $P$ and represented by
\begin{gather}\label{seis-1}
\text{$\ldots x_i, \ldots , x_3, x_2, x_1$}
\end{gather}
and
\begin{gather}\label{seis-2}
\text{$\ldots y_i, \ldots , y_3, y_2, y_1$}
\end{gather}
If there is a $j$ belonging to $N$ such that (i) $x_i = y_i$ for
$i> j$, (ii) $x_j$ is different from $y_j$ and (iii)
\begin{gather}\label{seis-3}
\text{$x_j, \ldots , x_3, x_2, x_1 > y_j, \ldots , y_3, y_2, y_1$}
\end{gather}
then the difference
\begin{gather}\label{seis-4}
\text{$(x_j, \ldots , x_3, x_2, x_1) - (y_j, \ldots , y_3, y_2,
y_1)$}
\end{gather}
is the finite distance between $x$ and $y$.
\end{definition}

\hspace{22pt} Let us notice that there is no way obtaining the sum
$x+y$ represented in the form (\ref{dois-10}) when $x$ and $y$
belong to $P$. For instance, be $x$ and $y$ identified by the same
structural-descriptive name

\begin{quote}
The represented number in the form (\ref{dois-10}) such that $f_i
= 0$ for $i$ odd and $f_i = 1$ for $i$ even.
\end{quote}

\hspace{22pt} Assuming that the sum $x+y$ is equal to
\begin{gather}\label{seis-5}
\text{$\ldots 1010101010100$}
\end{gather}
adding 1 the (\ref{seis-5})
\begin{gather}\label{seis-6}
\text{$\ldots 1010101010101$}
\end{gather}
and adding (\ref{seis-6}) to itself, we obtain
\begin{gather}\label{seis-7}
\text{$\ldots 0101010101010$}
\end{gather}
But (\ref{seis-7}) is $y$, therefore
\begin{gather}\label{seis-8}
\text{$((y+y)+1)+((y+y)+1)=y$}
\end{gather}

From (\ref{seis-8}) we concluded that it is false the hypothesis
that resulted in (\ref{seis-5}) and (\ref{seis-7}).

\section{The order of the elements}

\hspace{22pt} Let 0$\Phi$ be the sequence
\begin{gather}\label{sete-1}
\text{0$\mathfrak{SSS}\ldots\mathfrak{SSSS}SSSSSSSS\ldots$}
\end{gather}
such that
\begin{gather}\label{sete-2}
\begin{tabular}{c c c c c c c c c}
$\mathfrak{S}$ & $\mathfrak{S}$ & $\mathfrak{S}$ & $\mathfrak{S}$ & \ldots & $\mathfrak{S}$ & $\mathfrak{S}$ & $\mathfrak{S}$ & $\mathfrak{S}$ \\
$\updownarrow$ & $\updownarrow$ & $\updownarrow$ & $\updownarrow$ &   & $\updownarrow$ & $\updownarrow$ & $\updownarrow$ & $\updownarrow$\\
0 & 2 & 4 & 6 & \ldots & 7 & 5 & 3 & 1 \\
\end{tabular}
\end{gather}

shown in (\ref{dois-2}). We denote by $o_1$ the number represented
by 0$\mathfrak{SSS}\ldots\mathfrak{SSSS}$. Substituting $o_1$ in
(\ref{sete-1}), we have
\begin{gather}\label{sete-3}
\text{$o_{1} SSSSSSSS\ldots$}
\end{gather}

Let us indicate in the form $\mathfrak{S}$ a sequence of symbols S
on the right of $o_1$ in (\ref{sete-3}) such that the sequence of
S's in the form $\mathfrak{S}$ of (\ref{sete-4}) satisfies
(\ref{sete-2}).
\begin{gather}\label{sete-4}
\text{$o_{1} \mathfrak{SSS}\ldots\mathfrak{SSSS}SSSSSSSS\ldots$}
\end{gather}

Let us denote by $o_2$ the number represented by $o_{1}
\mathfrak{SSS}\ldots\mathfrak{SSSS}$. Repeating the process
obtains the succession
\begin{gather}\label{sete-5}
\text{$o_{1}, o_{2},o_{3},o_{4},\ldots$}
\end{gather}

We adopt the following axiom:
\begin{gather}\label{sete-6}
\textit{A6. The number w is the last element of P}
\end{gather}

\hspace{22pt} Let be
\begin{gather}\label{sete-7}
\text{$Q = P - \{ w\} $}
\end{gather}

Be $C_1$ the set of the elements of $Q$ that are at a
\textit{finite distance} to $o_1$. It is $C_2$ the set of the
elements of $Q$ that are to a finite distance of $o_2$. Be $C_n$
the set of the elements of $Q$ that are to a finite distance of
$o_n$, and so forth. With this ordination, we have
\begin{equation}\label{sete-8}
\begin{split}
M& =\{ 0,1,2,3,\ldots
,o_{1}^{-2},o_{1}^{-1},o_{1},o_{1}^{+1},o_{1}^{+2}, \ldots
,o_{2}^{-2},o_{2}^{-1},o_{2},o_{2}^{+1},o_{2}^{+2}, \ldots ,\\
& \quad o_{3}^{-2},o_{3}^{-1},o_{3},o_{3}^{+1},o_{3}^{+2}, \ldots
, w, w_1, w_3, w_3, \ldots \}
\end{split}
\end{equation}

Attributing to $o_i$ the exponent 0, we have
\begin{gather}\label{sete-9}
\text{$\{ \ldots
,o_{i}^{-2},o_{i}^{-1},o_{i}^{0},o_{i}^{+1},o_{i}^{+2}, \ldots \}
$}
\end{gather}

and abstracting the exponents of the $C_i$'s, we obtain
\begin{gather}\label{sete-10}
\text{$\{ \ldots ,-3,-2,-1,0,1,2,3, \ldots \} $}
\end{gather}

Let us notice that (\ref{sete-10}) is the set of the integer
numbers $Z$.

\section{Comment}

\hspace{22pt} In the literature $N$ denotes the set of the natural
numbers and to denominate $N$ as the set of the finite natural
numbers is a redundancy, because the existence of infinite natural
numbers is not considered. In this work, $N$ denotes the set of
the finite natural numbers and $M-N$ is the set of the infinite
natural numbers. Why not give the $M-N$ another denomination? Why
not use natural with some prefix to denominate the numbers of
$M-N$? With other denomination for $M-N$, we would avoid to
mention finite natural numbers when refering to the elements of
$N$, we would just say natural numbers . To reinforce the non
denomination of the elements of $M-N$ as natural numbers, there
are important properties that distinguish $M-N$ of $N$. For
instance, the elements of $N$ are archimedeans, that is,

\begin{quote}
If we took two numbers $m$ and $n$ belonging to $N$ such that $n>
m> 0$, it is always possible to add $m$ to itself a number $p$ of
times, where $p \in N$, such that the sum is greater than $n$.
\end{quote}

while the set $M-N$ does not possess this property. However, (i)
the most fundamental properties of $N$ (Peano's axioms) and of $M$
(Definition 2.2) are given by axioms of same form, (ii) so much
the elements of $N$ as the elements of $M-N$ are instances of $0X$
and (iii), for instance, the limit
\begin{equation}\label{oito-1}
\lim_{(n \in L, n \rightarrow \infty )}1 \ldots _{n} 11 = w
\end{equation}

is as valid as for non existence of the limit
\begin{equation}\label{oito-2}
\lim_{(n \in N, n \rightarrow \infty )}n
\end{equation}

\end{document}